 
\nopagenumbers
\baselineskip=14pt
\parskip=10pt
\def\Tilde{\char126\relax}
\def\halmos{\hbox{\vrule height0.15cm width0.01cm\vbox{\hrule height
 0.01cm width0.2cm \vskip0.15cm \hrule height 0.01cm width0.2cm}\vrule
 height0.15cm width 0.01cm}}
\font\eightrm=cmr8  
\font\eighttt=cmtt8
\magnification=\magstephalf

\parindent=0pt
\overfullrule=0in
 
\bf
\centerline
{
PROOF OF A CONJECTURE OF CHAN, ROBBINS, AND  YUEN
}
\rm
\bigskip
\centerline{ {\it 
Doron ZEILBERGER
}\footnote{$^1$}
{\eightrm  \raggedright
Department of Mathematics, Temple University,
Philadelphia, PA 19122, USA. 
{\eighttt zeilberg@math.temple.edu} \hfill \break
{\eighttt http://www.math.temple.edu/\Tilde zeilberg/   .}
Nov. 17, 1998.  Supported in part by the NSF.
} 
}
 
{\bf Abstract:}  Using the celebrated Morris Constant Term Identity,
we deduce a recent conjecture of Chan, Robbins, and Yuen
(math.CO/9810154), that asserts that the volume of a certain 
$n(n-1)/2$-dimensional polytope
is given in terms of the product of the first $n-1$ Catalan numbers.
 
Chan, Robbins, and Yuen[CRY]
conjectured that the cardinality of a certain set of triangular arrays
${\cal A}_n$ defined in pp. 6-7 of [CRY] equals the
product of the first $n-1$ Catalan numbers.
It is easy to see that their conjecture is equivalent
to the following {\it constant term identity}
(for any rational function $f(z)$ of a variable $z$,
$CT_{z} f(z)$ is the coeff. of $z^0$ in the formal Laurent
expansion of $f(z)$ (that always exists)):
$$
CT_{x_n} \dots CT_{x_1}
\prod_{i=1}^{n} (1-x_i)^{-2} \prod_{1 \leq i < j \leq n}
(x_j - x_i)^{-1}
= \prod_{i=1}^{n} {{1} \over {i+1}} {{2i} \choose {i}} \quad .
\eqno(CRY)
$$
But this is just the special case $a=2$, $b=0$, $c=1/2$, of the
{\it Morris Identity}[M] (where we made some trivial changes
of discrete variables, and `shadowed' it)
$$
CT_{x_n} \dots CT_{x_1}
\prod_{i=1}^{n} (1-x_i)^{-a} \prod_{i=1}^{n} x_i^{-b} 
\prod_{1 \leq i < j \leq n}
(x_j - x_i)^{-2c}
={{1} \over {n!}} \prod_{j=0}^{n-1}
{{\Gamma(a+b+(n-1+j)c)\Gamma(c)} \over
{\Gamma(a+jc)\Gamma(c+jc)\Gamma(b+jc+1)}}\,.
\eqno(Chip)
$$ 
To show that the right side of $(Chip)$ reduces to
the right side of $(CRY)$ upon the specialization
$a=2,b=0,c=1/2$, do the plugging in the former and call it
$M_n$. Then manipulate the products
to simplify $M_n/M_{n-1}$, and {\it then} use {\it Legendre's
duplication formula} 
$\Gamma(z) \Gamma(z+1/2)=\Gamma(2z)\Gamma(1/2)/2^{2z-1}$
three times, and {\it voil\`a}, up pops the Catalan
number  ${{2n} \choose {n}}/(n+1)$. \halmos
 
{\bf Remarks:} {\bf 1.} By converting the left side of $(Chip)$ into
a contour integral, we get the same integrand as in
the Selberg integral (with 
$a \rightarrow -a, \,b \rightarrow -b-1,\,c \rightarrow -c$). 
Aomoto's
proof of the Selberg
integral (SIAM J. Math. Anal. {\bf 18}(1987), 545-549) goes verbatim.
{\bf 2.} Conjecture 2 in [CRY] follows in the same way, from 
(the obvious contour-integral analog of) Aomoto's extension of 
Selberg's integral. Introduce a new
variable $t$, stick $CT_t t^{-k}$  in front of $(CRY)$, and
replace $(1-x_i)^{-2}$ by
$(1-x_i)^{-1}(t+x_i/(1-x_i))$. {\bf 3.} Conjecture 3 follows in the
same way from another specialization of $(Chip)$.
 
{\bf References}
\hfill\break
[CRY] Clara S. Chan, David P. Robbins, and David S. Yuen,
{\it On the volume of a certain polytope}, {\tt math.CO/9810154}.
\hfill\break
[M] Walter Morris, {\it ``Constant term identities for
finite and affine root systems, conjectures and theorems''},
Ph.D. thesis, University of Wisconsin, Madison, Wisconsin, 1982.
 
\bye